\documentclass{article}
\usepackage[dvips]{graphics}
\usepackage[english]{babel}
\usepackage{amsfonts}
\usepackage{amsthm}
\usepackage{amsmath}
\usepackage{amscd}
\usepackage[active]{srcltx} % SRC Specials: DVI [Inverse] Search
\usepackage{latexsym}
\usepackage{amssymb}
\usepackage[latin1]{inputenc}

\newcommand{\matriz}[4]{\displaystyle\
    \left(
       \begin{array}{cc}
        {#1}&{#2}\\
        {#3}&{#4}
       \end{array}
     \right)}

\newcommand{\PI}[2]{\left\langle \,#1 , #2\, \right\rangle}
\newcommand{\K}[2]{[ \,#1 , #2\, ]}

\newcommand{\ra}{\rightarrow}

\newcommand{\x}{\times}
\newcommand{\s}{\sigma}

\newcommand{\CC}{\mathbb{C}}

\newcommand{\NN}{\mathbb{N}}
\newcommand{\St}{\mathcal{S}}
\newcommand{\Q}{\mathcal{Q}}
\newcommand{\T}{\mathcal{T}}
\newcommand{\M}{\mathcal{M}}
\newcommand{\N}{\mathcal{N}}
\newcommand{\U}{\mathcal{U}}
\newcommand{\HH}{\mathcal{H}}
\newcommand{\KK}{\mathcal{K}}
\newcommand{\mc}[1]{\mathcal{#1}}
\newcommand{\F}{\mathcal{F}}

\newcommand{\ol}{\overline}

\newcommand{\ort}{[\bot]}

\newcommand{\vp}{\varphi}

\newcommand{\noi}{\noindent}

\DeclareMathOperator{\Span}{span} 
\DeclareMathOperator{\ind}{ind} 
\DeclareMathOperator{\sgn}{sgn}

\newtheorem{thm}{Theorem}[section]

\newtheorem{prop}[thm]{Proposition}

\newtheorem{lem}[thm]{Lemma}

\newtheorem{cor}[thm]{Corollary}

\theoremstyle{definition}

\newtheorem{rem}[thm]{Remark}

\newtheorem*{pf*}{Proof}

\newtheorem*{defn}{Definition}

\newtheorem{exmp}{Example}

\begin{document}

\title{On a family of frames for Krein spaces}

\author{J. I. Giribet, A. Maestripieri, F. Mart\'{\i}nez Per\'{\i}a and P. Massey}

\date{}

\maketitle

\begin{abstract}
A definition of frames for Krein spaces is proposed, which extends the
notion of $J$-orthonormal basis of Krein spaces. A $J$-frame for a Krein space $(\HH, \K{\,}{\,})$ is in particular a frame for $\HH$ in the Hilbert space sense. But it is also compatible with the indefinite inner product $\K{\,}{\,}$, meaning that it determines a pair of maximal uniformly $J$-definite subspaces with different positivity, an analogue to the maximal dual pair associated to a $J$-orthonormal basis. 

Also, each $J$-frame induces an indefinite reconstruction formula for the vectors in $\HH$, which resembles the one given by a $J$-orthonormal basis.
\end{abstract}

keywords: Krein spaces, \ frames, \ uniformly $J$-definite subspaces


\section{Introduction}

In recent years, frame theory for Hilbert spaces has been thoroughly developed, see e. g. \cite{TaF, Chr, Chr2, HanLarson}. Fixed a Hilbert space $(\HH,\PI{\,}{\,})$, a frame for $\HH$ is a (generally overcomplete) family of vectors $\F=\{f_i\}_{i\in I}$ in $\HH$ which satisfies the inequalities
\begin{equation}\label{frame}
	A\|f\|^2\leq \sum_{i\in I} |\PI{f}{f_i}|^2 \leq B\|f\|^2, \ \ \ \text{for every $f\in \HH$},
\end{equation}
for positive constants $0<A\leq B$. The (bounded, linear) operator $S:\HH\ra\HH$ defined by
\begin{equation}\label{frameop}
Sf=\sum_{i\in I}\PI{f}{f_i}f_i, \ \ \ f\in \HH,
\end{equation}	
is known as the frame operator associated to $\F$. The inequalities in Eq. \eqref{frame} imply that $S$ is a (positive) boundedly invertible operator, and it allows to reconstruct each vector $f\in\HH$ in terms of the family $\F$ as follows: 
\begin{equation}
	f=\sum_{i\in I}\PI{f}{S^{-1}f_i}f_i=\sum_{i\in I}\PI{f}{f_i}S^{-1}f_i.
\end{equation}
The above formula is known as the \emph{reconstruction formula associated to $\F$}. Notice that if $\F$ is a Parseval frame, i.e. if $S=I$, then the reconstruction formula resembles the Fourier series of $f$ associated to an orthonormal basis $\mc{B}=\{b_k\}_{k\in K}$ of $\HH$:
\[
f=\sum_{k\in K}\PI{f}{b_k}b_k,
\]
but the \emph{frame coefficients} $\{\PI{f}{f_i}\}_{i\in I}$ given by $\F$ allow to reconstruct $f$ even when some of these coefficients are missing (or corrupted). Indeed, each vector $f\in\HH$ may admit several reconstructions in terms of the frame coefficients as a consequence of the redundancy of $\F$. These are some of the advantages of frames over (orthonormal, orthogonal or Riesz) bases in signal processing applications, when noisy channels are involved, e.g. see \cite{BodPau,HolPau,Stro}.

Given a Krein space $(\HH,\K{\,}{\,})$ with fundamental symmetry $J$, a \emph{$J$-orthonormalized system} is a family $\mc{E}=\{e_i\}_{i\in I}$ such that $\K{e_i}{e_j}=\pm\delta_{ij}$, for $i, j\in I$. A $J$-orthonormal basis is a $J$-orthonormalized system which is also a Schauder basis for $\HH$. If $\mc{E}=\{e_i\}_{i\in I}$ is a $J$-orthonormal basis of $\HH$ then the vectors in $\HH$ can be represented as follows:
\begin{equation}
	f=\sum_{i\in I} \s_i\ \K{f}{e_i}\ e_i, \ \ \ f\in\HH,
\end{equation}
where $\s_i=\K{e_i}{e_i}=\pm 1$.

$J$-orthonormalized systems (and bases) are intimately related to the notion of dual pair. In fact, each $J$-orthonormalized system generates a \emph{dual pair}, i.e. a pair $(\mc{L}_+,\mc{L}_-)$ of subspaces of $\HH$ such that $\mc{L}_+$ is $J$-nonnegative, $\mc{L}_-$ is $J$-nonpositive and $\mc{L}_+$ is $J$-orthogonal to $\mc{L}_-$, i.e. $\K{\mc{L}_+}{\mc{L}_-}=0$. Moreover, if $\mc{E}$ is a $J$-orthonormal basis of $\HH$, the dual pair associated to $\mc{E}$ is maximal (with respect to the inclusion preorder) and the subspaces $\mc{L}_+$ and $\mc{L}_-$ are uniformly $J$-definite, see \cite[Ch.1, \S 10]{Iokhvidov}. Therefore the dual pair $(\mc{L}_+,\mc{L}_-)$ is a fundamental decomposition of $\HH$. Notice that, considering the Hilbert space structure induced by the above fundamental decomposition, the $J$-orthonormal basis $\mc{E}$ turns out to be an orthonormal basis in the associated Hilbert space. Therefore, each $J$-orthonormal basis can be realized as an orthonormal basis of $\HH$ (respect to an appropriate definite inner product).

Given a pair of maximal uniformly $J$-definite subspaces $\M_+$ and $\M_-$ of a Krein space $\HH$, with different positivity, if $\F_\pm=\{f_i\}_{i\in I_\pm}$ is a frame for the Hilbert space $(\M_\pm, \pm\K{\,}{\,})$, it is easy to see that 
\[
\F=\F_+\cup\F_-,
\]
is a frame for $\HH$, which produces an \emph{indefinite reconstruction formula}:
\[
f= \sum_{i\in I} \s_i \K{f}{g_i}f_i=\sum_{i\in I} \s_i \K{f}{f_i}g_i, \ \ \ f\in\HH,
\]
where $\s_i=\sgn\K{f_i}{f_i}$ and $\{g_i\}_{i\in I}$ is some (equivalent) frame for $\HH$ (see Example \ref{glueing} and Proposition \ref{J dual}).

\medskip

The aim of this work is to introduce and characterize a particular family of frames for a Krein space $(\HH,\K{\,}{\,})$ -hereafter called $J$-frames- that are compatible with the indefinite inner product $\K{\,}{\,}$. Some different approaches to frames for Krein spaces and indefinite reconstruction formulas are developed in \cite{Wagner} and \cite{Waldron}, respectively.

\medskip

The paper is organized as follows: Section 2 contains some
preliminaries results both in Krein spaces and in frame theory for Hilbert spaces. 

Section 3 presents the $J$-frames. Briefly,  a $J$-frame for the Krein space $(\HH,\K{\,}{\,})$ is a Bessel family $\F=\{f_i\}_{i\in I}$ with synthesis operator $T:\ell_2(I)\ra\HH$ such that the ranges of $T_+:=TP_+$ and $T_-:=T(I-P_+)$ are maximal uniformly $J$-positive and maximal uniformly $J$-negative subspaces, respectively, where $I_+=\{i\in I: \K{f_i}{f_i}>0 \}$ and $P_+$ is the orthogonal projection onto $\ell_2(I_+)$, as a subspace of $\ell_2(I)$.

It is immediate that $J$-orthonormal bases are $J$-frames, because they generate maximal dual pairs \cite[Ch. 1, \S 10.12]{Iokhvidov}.

Also, if $\F$ is a $J$-frame for $\HH$, observe that $R(T)=R(T_+)+ R(T_-)$ and recall that the sum of a pair of maximal uniformly $J$-definite subspaces with different positivity coincides with $\HH$ \cite[Corollary 1.5.2]{Ando}. Therefore, each $J$-frame is in fact a frame for $\HH$ in the Hilbert space sense. Moreover, it is shown that $\F_+=\{f_i\}_{i\in I_+}$ is a frame for the Hilbert space $(R(T_+),\K{\,}{\,})$ and $\F_-=\{f_i\}_{i\in I\setminus I_+}$ is a frame for $(R(T_-),-\K{\,}{\,})$, i.e. there exist constants $B_-\leq A_-<0<A_+\leq B_+$ such that
	\begin{equation}
	A_\pm \K{f}{f} \leq \sum_{i\in I_\pm} |\K{f}{f_i}|^2 \leq B_\pm \K{f}{f} \ \ \ \text{for every $f\in R(T_\pm)$}.
	\end{equation}
The optimal constants satisfying the above inequalities can be characterized in terms of $T_\pm$ and the Gramian operators of their ranges.

This section ends with a geometrical characterization of $J$-frames, in terms of the (minimal) angles between the uniformly $J$-definite subspace $R(T_\pm)$ and the cone of neutral vectors of the Krein space.

Section 4 is devoted to study the synthesis operators associated to $J$-frames. Fixed a Krein space $\HH$ and given a bounded operator $T:\ell_2(I)\ra\HH$, it is described under which conditions $T$ is the synthesis operator of a $J$-frame.

In Section 5 the $J$-frame operator is introduced. Given a $J$-frame $\mc{F}=\{f_i\}_{i\in I}$, the \emph{$J$-frame operator} $S:\HH\ra \HH$ is defined by
\[
Sf=\sum_{i\in I}\s_i \ \K{f}{f_i} \ f_i, \ \ \ f\in \HH,
\]
where $\s_i=\sgn(\K{f_i}{f_i})$. This operator resembles the frame operator for frames in Hilbert spaces (see Eq. \eqref{frameop}), and it has similar properties, in particular $S=TT^\#$ if $T:\ell_2(I)\ra\HH$ is the synthesis operator of $\F$ (see Proposition \ref{props S}). Furthermore, each $J$-frame $\F=\{f_i\}_{i\in I}$ determines an \emph{indefinite reconstruction formula}, which depends on the $J$-frame operator $S$:
\begin{equation}
	f= \sum_{i\in I} \s_i\ \K{f}{S^{-1}f_i}\ f_i = \sum_{i\in I} \s_i\ \K{f}{f_i}\ S^{-1}f_i, \ \ \ \text{for every $f\in \HH$}.
\end{equation}
In this case the family $\{S^{-1}f_i\}_{i\in I}$ turns out to be a $J$-frame too.

Finally, it will be shown that the $J$-frame operator of a $J$-frame $\F$ is intimately related to the projection $Q=P_{R(T_+)//R(T_-)}$ determined by the decomposition $\HH=R(T_+)\dotplus R(T_-)$. In fact, fixed a $J$-selfadjoint invertible operator $S$ acting on a Krein space $\HH$, it is the $J$-frame operator for a $J$-frame $\F$ if and only if there exists a projection $Q$ with uniformly $J$-definite range and kernel such that $QS$ is a $J$-positive operator and $(I-Q)S$ is a $J$-negative operator, see Theorem \ref{carac S}.

\section{Preliminaries}

Along this work $\HH$ denotes a complex (separable) Hilbert space. If $\mc{K}$ is another Hilbert space then  $L(\HH, \mc{K})$ is the algebra of bounded linear operators from $\HH$ into $\mc{K}$ and $L(\HH)=L(\HH,\HH)$. The groups of linear invertible and unitary operators acting on $\HH$ are denoted by $GL(\HH)$ and $\U(\HH)$, respectively. Also,  $L(\HH)^+$ denotes the cone of positive semidefinite operators acting on $\HH$ and $GL(\HH)^+=GL(\HH)\cap L(\HH)^+$.

If $T\in L(\HH, \mc{K})$ then $T^*\in L(\mc{K},\HH)$ denotes the adjoint operator of $T$, $R(T)$ stands for its range and $N(T)$ for its nullspace. Also, if $T\in L(\HH, \mc{K})$ has closed range, $T^\dagger\in L(\KK,\HH)$ denotes the Moore-Penrose inverse of $T$.

Hereafter, $\St\dotplus \T$ denotes the direct sum of two (closed) subspaces $\St$ and $\T$ of $\HH$. On the other hand, $\St\oplus \T$ stands for the (direct) orthogonal sum of them and $\St\ominus \T:= \St\cap (\St\cap\T)^\bot$. If $\HH=\St\dotplus \T$, the oblique projection onto $\St$ along $\T$ is the unique projection with range $\St$ and nullspace $\T$. It is denoted by $P_{\St//\T}$. In particular, $P_\St:=P_{\St// \St^\bot}$ is the orthogonal projection onto $\St$.

\subsection{Krein spaces}

In what follows we present the standard notation and some basic results on Krein spaces. For a complete exposition on the subject (and the proofs of the results below) see the books by J. Bogn\'ar \cite{Bognar} and T. Ya. Azizov and I. S. Iokhvidov \cite{Iokhvidov} and the monographs by T. Ando \cite{Ando} and by M. Dritschel and J. Rovnyak \cite{Dritschel 1}.

Given a Krein space $(\HH, \K{\,}{\,})$ with a \emph{fundamental decomposition} $\HH=\HH_+\dotplus \HH_-$, the direct (orthogonal) sum of the Hilbert spaces $(\HH_+, \K{\,}{\,})$ and $(\HH_-, -\K{\,}{\,})$ is denoted by $(\HH,\PI{\,}{\,})$. 

Observe that the indefinite metric and the inner product of $\HH$ are related by means of a  \emph{fundamental symmetry}, i.e. a unitary selfadjoint operator $J\in L(\HH)$ which satisfies:
\[
\K{x}{y}=\PI{Jx}{y}, \ \ \ \ x,y\in\HH.
\]
If $\HH$ and $\KK$ are Krein spaces, $L(\HH,\KK)$ stands for the vector space of linear transformations which are bounded respect to the associated Hilbert spaces $(\HH, \PI{\,}{\,}_\HH)$ and $(\KK, \PI{\,}{\,}_\KK)$. Given $T\in L(\HH,\KK)$, the $J$-adjoint operator of $T$ is defined by $T^\#=J_\HH T^* J_\KK$, where $J_\HH$ and $J_\KK$ are the fundamental symmetries associated to $\HH$ and $\KK$, respectively. An operator $T\in L(\HH)$ is $J$-selfadjoint if $T=T^\#$.

A vector $x\in \HH$ is \emph{$J$-positive} if $\K{x}{x}>0$. A subspace $\St$ of $\HH$ is \emph{$J$-positive} if every $x\in\St$, $x\neq 0$, is a $J$-positive vector. A subspace $\St$ of $\HH$ is \emph{uniformly $J$-positive} if there exists $\alpha> 0$ such that
\[
\K{x}{x} \geq \alpha \|x\|^2, \ \ \ \ \text{for every $x\in\St$},
\]
where $\|\, \|$ stands for the norm of the associated Hilbert space $(\HH,\PI{\,}{\,})$. 

$J$-nonnegative, $J$-neutral, $J$-negative, $J$-nonpositive and uniformly $J$-negative vectors and subspaces are defined analogously.  

\begin{rem}\label{Pedro}
If $\St_+$ is a closed uniformly $J$-positive subspace of a Krein space $(\HH,\K{\,}{\,})$, observe that $(\St_+, \K{\,}{\,})$ is a Hilbert space. In fact, the forms $\K{\,}{\,}$ and $\PI{\,}{\,}$ are equivalent inner products on $\St_+$, because 
\[
\alpha \|f\|^2\leq \K{f}{f} \leq \|f\|^2, \ \ \ \text{for every $f\in\St_+$}.
\]
Analogously, if $\St_-$ is a closed uniformly $J$-negative subspace of $(\HH,\K{\,}{\,})$, $(\St_-, -\K{\,}{\,})$ is a Hilbert space.
\end{rem}

\begin{prop}[\cite{Iokhvidov}, Cor. 7.17]\label{LSS}
Let $\HH$ be a Krein space with fundamental symmetry $J$ and $\St$ a $J$-nonnegative closed subspace of $\HH$. Then, $\St$ is the range of a $J$-selfadjoint projection if and only if $\St$ is uniformly $J$-positive.
\end{prop}

\noi Recall that, given a closed subspace $\M$ of a Krein space $\HH$, the Gramian operator of $\M$ is defined by:
\[
G_\M=P_\M J P_\M,
\]
where $P_\M$ is the orthogonal projection onto $\M$ and $J$ is the fundamental symmetry of $\HH$. If $\M$ is  $J$-semidefinite, then $\M\cap\M^{\ort}$ coincides with $\N:=\{f\in \M: \K{f}{f}=0\}$. Therefore, it is easy to see that
\[
G_\M=G_{\M\ominus\N}.
\]
Given a subspace $\St$ of a Krein space $\HH$, the \emph{$J$-orthogonal companion} to $\St$ is defined by 
\[
\St^{\ort}=\{ x\in\HH : \K{x}{s}=0 \; \text{for every $s\in\St$}\}.
\]
A subspace $\St$ of $\HH$ is \emph{$J$-non degenerated} if $\St\cap\St^{\ort}=\{0\}$. Notice that if $\St$ is a $J$-definite subspace of $\HH$ then it is $J$-non degenerated.

\subsection{Angles between subspaces and reduced minimum modulus}

Given two closed subspaces $\St$ and $\T$ of a Hilbert space $\HH$, the cosine of the \emph{Friedrichs angle} between $\St$ and $\T$ is defined by
\[
c(\St,\T)=\sup\{ |\PI{x}{y}| : x\in\St\ominus\T, \|x\|=1, \, y\in\T\ominus\St, \|y\|=1\}.
\]

\noi It is well known that
\[
c(\St,\T)<1 \ \ \Leftrightarrow \ \ \St + \T \ \text{is closed} \ \ \Leftrightarrow \ \ c(\St^\bot,\T^\bot)<1.
\]
Furthermore, if $P_\St$ and $P_\T$ are the orthogonal projections onto $\St$ and $\T$, respectively, then $c(\St,\T)<1$ if and only if $(I-P_\St)P_\T$ has closed range.
See \cite{Deu} for further details.

\medskip

\noi The next definition is due to T. Kato, see \cite[Ch. IV, $\S$ 5]{Kato}.	
\begin{defn}
The \emph{reduced minimum modulus} $\gamma(T)$ of an operator $T\in L(\HH,\KK)$ is defined
by
\[
\gamma(T) = \inf\{\|Tx\| : \ x \in N(T)^\bot, \, \|x\| = 1\}.
\]
\end{defn}
Observe that $\gamma(T)=\sup\{C\geq 0 : C \|x\| \leq \|Tx\| \ \text{for every} \ x \in N(T)^\bot, \, \|x\| = 1\}$. It is well known that $\gamma(T) = \gamma(T^*) = \gamma(T^*T)^{1/2}$. Also, it can be shown that an operator $T\neq 0$ has closed range if and only if $\gamma(T) > 0$. In this case, $\gamma(T) = \|T^\dag\|^{-1}$.

If $\HH$ and $\KK$ are Krein spaces with fundamental symmetries $J_\HH$ and $J_\KK$, respectively, and $T\in L(\HH,\KK)$ then
\[
\gamma(T^\#)=\gamma(J_\HH T^* J_\KK)=\gamma(T^*)=\gamma(T),
\]
because $J_\HH$ (resp. $J_\KK$) is a unitary operator on $\HH$ (resp. $\KK$). 

\begin{rem}\label{prop gama}
If $\M_+$ is a closed $J$-nonnegative subspace of a Krein space $\HH$ then
\begin{equation}
\gamma(G_{\M_+})=\alpha^+,
\end{equation}
where $\alpha^+\in [0,1]$ is the supremum among the constants $\alpha\in [0,1]$ such that $\alpha \|f\|^2 \leq \K{f}{f}$ for every $f\in\M_+$. From now on, the constant $\alpha^+$ is called the \emph{definiteness bound} of $\M_+$. Notice that $\alpha^+$ is in fact a maximum for the above set and $\M^+$ is uniformly $J$-positive if and only if $\alpha^+>0$. 

Analogously, if $\M_-$ is a  $J$-nonpositive subspace then $\gamma(G_{\M_-})=\alpha^-$, where $\alpha^-$ is the definiteness bound of $\M_-$, i.e.
\[ 
\alpha^- =\max\{\alpha\in[0,1]:\ \K{f}{f}\leq -\alpha\,\|f\|^2 \ \ \text{for every $f\in \M_-$} \}.
\]
\end{rem}

\subsection{Frames for Hilbert spaces}

\noi The following is the standard notation and some basic results on frames for Hilbert spaces, see \cite{TaF, Chr, HanLarson}.

\medskip

A \emph{frame} for a Hilbert space $\HH$ is a family of vectors $\mc F=\{f_i\}_{i\in I}\subset \HH$ for which there exist constants $0<A\leq B<\infty$ such that 
\begin{equation}\label{ecu frames} 
A\ \|f\|^2 \leq \sum_{i\in I} |\langle f,f_i\rangle |^2 \leq B\ \|f\|^2\, , \ \ \text{for every $f\in \mc{H}$}.
\end{equation}
The optimal constants (maximal for $A$ and minimal for $B$) are known, respectively, as the upper and lower frame bounds. 

If a family of vectors $\mc F=\{f_i\}_{i\in I}$ satisfies the upper bound condition in \eqref{ecu frames}, then $\mc{F}$ is a \emph{Bessel family}. For a Bessel family $\mc F=\{f_i\}_{i\in I}$, the \emph{synthesis operator} $T\in L(\ell_2(I),\HH)$ is defined by
\[
Tx=\sum_{i\in I}\PI{x}{e_i}f_i,
\] 
where $\{e_i\}_{i\in I}$ is the standard orthonormal basis of $\ell_2(I)$. It holds that $\mc{F}$ is a frame for $\HH$ if and only if $T$ is surjective. In this case, the operator $S=TT^*\in L(\HH)$ is invertible and is called the \emph{frame operator}. It can be easily verified that
\begin{equation}\label{ecu S}
Sf=\sum_{i\in I}\PI{f}{f_i}f_i, \ \ \ \text{for every $f\in\HH$}.
\end{equation}
This implies that the frame bounds can be computed as: $A=\|S^{-1}\|^{-1}$ and $B=\|S\|$. From \eqref{ecu S}, it is also easy to obtain the \emph{canonical reconstruction formula} for the vectors in $\HH$:
\[
f=\sum_{i\in I}\PI{f}{S^{-1}f_i}f_i= \sum_{i\in I}\PI{f}{f_i}S^{-1}f_i, \ \ \ \text{for every $f\in\HH$},
\] 
and the frame $\{S^{-1}f_i\}_{i\in I}$ is called the \emph{canonical dual frame} of $\mc{F}$. More generally, if a frame $\mc{G}=\{g_i\}_{i\in I}$ satisfies 
\begin{equation}\label{dual}
f=\sum_{i\in I}\PI{f}{g_i}f_i= \sum_{i\in I}\PI{f}{f_i}g_i, \ \ \ \text{for every $f\in\HH$},
\end{equation}
then $\mc{G}$ is called a \emph{dual frame} of $\mc{F}$.

\section{$J$-frames: definition and basic properties}

Let $\HH$ be a Krein space with fundamental symmetry $J$. Given a Bessel family $\mc{F}=\{f_i\}_{i\in I}$ in $\HH$ consider the synthesis operator $T\in L(\ell_2(I),\HH)$. If $I_+=\{i\in I: \ \K{f_i}{f_i}\geq 0\}$ and $I_-=\{i\in I: \ \K{f_i}{f_i}< 0\}$, consider the orthogonal decomposition of $\ell_2(I)$ given by
\begin{equation}\label{desc fund}
\ell_2(I)=\ell_2(I_+)\oplus \ell_2(I_-),
\end{equation}
and denote by $P_\pm$ the orthogonal projection onto $\ell_2(I_\pm)$. Also, let $T_\pm=TP_\pm$. If $\M_\pm=\ol{\Span\{ f_i : \ i\in I_\pm\}}$, notice that $\Span\{ f_i : \ i\in I_\pm\}\subseteq R(T_\pm) \subseteq \M_\pm$ and 
\[
R(T)=R(T_+) + R(T_-).
\]

\begin{defn}\label{Jframe}
The Bessel family $\mc{F}=\{f_i\}_{i\in I}$ is a \emph{$J$-frame} for $\HH$ if $R(T_+)$ is a maximal uniformly $J$-positive subspace of $\HH$ and $R(T_-)$ is a maximal uniformly $J$-negative subspace of $\HH$.
\end{defn}

\medskip

Notice that, in particular, every $J$-orthogonalized basis of a Krein space $\HH$ is a $J$-frame for $\HH$, because it generates a maximal dual pair, see \cite[Ch. 1, \S 10.12]{Iokhvidov}. 

If $\F$ is a $J$-frame, as a consequence of its maximality, $R(T_\pm)$ is closed. So, $R(T_\pm)=\M_\pm$ and, by \cite[Corollary 1.5.2]{Ando}, $\M_+ + \M_-=\HH$. Then, it follows that $\mc{F}$ is a frame for the associated Hilbert space $(\HH,\PI{\,}{\,})$ because 
\[
R(T)=R(T_+) + R(T_-)=\M_+ + \M_-=\HH.
\]

Given a Bessel family $\F=\{f_i\}_{i\in I}$, consider the subspaces $R(T_+)$ and $R(T_-)$ as above. If $K_\pm: \mc{D_\pm}\ra \HH_\mp$ is the angular operator associated to $R(T_\pm)$, the \emph{operator of transition} associated to the Bessel family $\F$ is defined by
\[
F=K_+P + K_-(I-P) : \mc{D}_+ + \mc{D}_- \ra \HH,
\]
where $P=\frac{1}{2}(I+J)$ is the $J$-selfadjoint projection onto $\HH_+$ and $\mc{D}_\pm$ is a subspace of $\HH_\pm$ (the domain of $K_\pm$), see \cite{Opus}. 
\begin{prop}
Let  $\F=\{f_i\}_{i\in I}$ be a Bessel family in $\HH$. Then, $\F$ is a $J$-frame if and only if $F$ is everywhere defined (i.e. $\mc{D}_+ + \mc{D}_-=\HH$) and $\|F\|<1$.
\end{prop}

\begin{pf*}{Proof}{}
See \cite[Proposition 2.6]{Opus}.
  \end{pf*}

It follows from the definition that, given a $J$-frame $\mc{F}=\{f_i\}_{i\in I}$ for the Krein space $\HH$, $\K{f_i}{f_i}\neq 0$ for every $i\in I$, i.e. $I_\pm=\{i\in I: \ \pm\K{f_i}{f_i}> 0\}$. This fact allows to endow the coefficients space $\ell_2(I)$ with a Krein space structure. Denote $\s_i= \sgn(\K{f_i}{f_i})=\pm 1$ for every $i\in I$. Then, the diagonal operator $J_2\in L(\ell_2(I))$ defined by
\begin{equation}\label{J2}
J_2\, e_i= \s_i\, e_i, \ \ \ \text{for every $i\in I$},
\end{equation}
is a selfadjoint involution on $\ell_2(I)$. Therefore, $\ell_2(I)$ with the fundamental symmetry $J_2$ is a Krein space.

Now, if $T\in L(\ell_2(I), \HH)$ is the synthesis operator of $\mc{F}$, the $J$-adjoints of $T$, $T_+$ and $T_-$ can be easily calculated, in fact if $f\in \HH$:
\[
T_\pm^\# f = \pm \sum_{i\in I_\pm} \K{f}{f_i}e_i,
\]
and $T^{\#}f=(T_+ + T_-)^\# f= T_+^{\#}f + T_-^{\#}f= \sum_{i\in I_+} \K{f}{f_i}e_i - \sum_{i\in I_-} \K{f}{f_i}e_i=\sum_{i\in I}\s_i\K{f}{f_i}e_i$.

\medskip

\begin{exmp}\label{ejem 1}
It is easy to see that not every frame of $J$-nonneutral vectors is a $J$-frame: given the Krein space obtained by endowing $\CC^3$ with the sesquilinear form 
\[
[(x_1,x_2,x_3),(y_1,y_2,y_3) ]= x_1\overline{y_1} + x_2\overline{y_2}  -x_3\overline{y_3},
\]
consider $f_1=(1,0,\frac{1}{\sqrt{2}})$, $f_2=(0,1,\frac{1}{\sqrt{2}})$ and $f_3=(0, 0,1)$. Observe that $\mc{F}=\{f_1,f_2,f_3\}$ is a frame for $\CC^3$ because it is a (linear) basis for the space.

On the other hand, $\M_+=\Span \{f_1,\,f_2\}$ and $\M_-=\Span\{f_3\}$. If $(a,b,\frac{1}{\sqrt{2}}(a+b))$ is an arbitrary vector in $\M_+$ then 
\[
\K{f}{f}=|a|^2+|b|^2-\tfrac{1}{2} |a+b|^2 =\tfrac{1}{2}|a-b|^2\geq 0,
\] 
so $\M_+$ is a $J$-nonnegative subspace of $\CC^3$. But $\M_+$ is not uniformly $J$-positive, because $(1,1,\sqrt{2})\in\M_+$ is a (non trivial) $J$-neutral vector. Therefore, $\mc{F}$ is not a $J$-frame for $(\CC^3,\K{\,}{\,})$. 
\end{exmp}

\medskip

The following is a handy way to construct $J$-frames for a given Krein space. Along this section, it will be shown that every $J$-frame can be realized in this way.

\begin{exmp}\label{glueing}
Given a Krein space $\HH$ with fundamental symmetry $J$, let $\M_+$ (resp. $\M_-$) be a maximal uniformly $J$-positive (resp. $J$-negative) subspace of $\HH$. If $\mc{F}_\pm=\{f_i\}_{i\in I_\pm}$ is a frame for the Hilbert space $(\M_\pm, \pm\K{\,}{\,})$ then $\mc{F}=\mc{F}_+ \cup \mc{F}_-$ is a $J$-frame for $\HH$.

Indeed, by Remark \ref{Pedro}, $\F_+$ and $\F_-$ are Bessel families in $\HH$. Hence, $\F$ is a Bessel family and, if $I=I_+ \dot{\cup} I_-$ (the disjoint union of $I_+$ and $I_-$), the synthesis operator $T\in L(\ell_2(I), \HH)$ of $\mc{F}$ is given by
\[
Tx = T_+ x_+ + T_- x_- \ \ \ \ \text{if} \ \ \ \ x=x_+ + x_-\in \ell_2(I_+)\oplus\ell_2(I_-)=:\ell_2(I),
\]
where $T_\pm: \ell_2(I_\pm)\ra \M_\pm$ is the synthesis operator of $\mc{F}_\pm$. Then, it is clear that $R(TP_\pm)=\M_\pm$ is a maximal uniformly $J$-definite subspace of $\HH$.
\end{exmp}

\medskip

\begin{prop}\label{jframes -> frames ineq}
Let $\mc{F}=\{f_i\}_{i\in I}$ be a $J$-frame for $\HH$. Then, $\mc{F}_\pm=\{f_i\}_{i\in I_\pm}$ is a frame for the Hilbert space $(\M_\pm, \pm\K{\,}{\,})$, i.e. there exist constants $B_-\leq A_-<0<A_+\leq B_+$ such that
	\begin{equation}\label{frm pos}
	A_\pm \K{f}{f} \leq \sum_{i\in I_\pm} |\K{f}{f_i}|^2 \leq B_\pm \K{f}{f} \ \ \ \text{for every $f\in \M_\pm$}.
	\end{equation}
\end{prop}

\begin{pf*}{Proof}{}
If  $\mc{F}=\{f_i\}_{i\in I}$ is a $J$-frame for $\HH$, then $R(T_+)=\M_+$ is a (maximal) uniformly $J$-positive subspace of $\HH$. So, $T_+$ is a surjection from $\ell_2(I)$ onto the Hilbert space $(\M_+,\K{\,}{\,})$. Therefore, $\F_+$ is a frame for $(\M_+,\K{\,}{\,})$. In particular, there exist constants $0<A_+\leq B_+$ such that Eq. \eqref{frm pos} is satisfied for $\M_+$. The assertion on $\F_-$ follows analogously.
  \end{pf*}

Now, assuming that $\F$ is a $J$-frame for a Krein space $(\HH,\K{\,}{\,})$, a set of constants $\{B_-, A_-, A_+, B_+\}$ satisfying Eq. \eqref{frm pos} is going to be computed. They depend only on the definiteness bounds for $R(T_\pm)$, the norm and the reduced minimum modulus of $T_\pm$.

Suppose that $\F$ is a $J$-frame for a Krein space $(\HH,\K{\,}{\,})$  with synthesis operator $T\in L(\ell_2(I),\HH)$. Since $R(T_+)=\M_+$ is a (maximal) uniformly $J$-positive subspace of $\HH$, there exists $\alpha_+>0$ such that $\alpha_+ \|f\|^2 \leq \K{f}{f}$ for every $f\in\M_+$. So,
\[
 \sum_{i\in I_+} |\K{f}{f_i}|^2 = \|T_+^\# f\|^2 \leq \|T_+^\#\|^2 \|f\|^2 \leq B_+ \K{f}{f}, \ \ \ \text{ for every $f\in \M_+$},
\]
where $B_+=\frac{\|T_+^\#\|^2}{\alpha_+}=\frac{\|T_+\|^2}{\alpha_+}$. Furthermore, since $N(T_+^\#)^\bot=J(\M_+)$, if $f\in\M_+$,
\begin{eqnarray*}
\sum_{i\in I_+} |\K{f}{f_i}|^2 &=& \|T_+^\# f\|^2= \|T_+^\# P_{J(\M_+)}f\|^2 \geq  \gamma(T_+^\#)^2 \|P_{J(\M_+)}f\|^2 = \gamma(T_+)^2 \|P_{\M_+}Jf\|^2=  \\ &=& \gamma(T_+)^2 \|G_{\M_+}f\|^2 \geq \gamma(T_+)^2 \gamma(G_{\M_+})^2 \|f\|^2 \geq A_+ \K{f}{f}, 
\end{eqnarray*}
where $A_+= \gamma(T_+)^2\gamma(G_{\M_+})^2= \gamma(T_+)^2\alpha_+^2$, see Remark \ref{prop gama}.

Analogously, $A_-=-\gamma(T_-)^2\alpha_-^2$ and $B_-=-\frac{\|T_-\|^2}{\alpha_-}$ satisfy Eq \eqref{frm pos} for every $f\in R(T_-)=\M_-$, if $\alpha_-$ is the definiteness bound of the (maximal) uniformly $J$-negative subspace $\M_-$. 

Usually, the bounds $A_\pm=\pm\alpha_\pm^2\gamma(T_\pm)^2$ and $B_\pm=\pm \frac{\|T_\pm\|^2}{\alpha_\pm}$ are not optimal for the $J$-frame $\F$.

\begin{defn} 
Let $\mc{F}$ be a $J$-frame for the Krein space $\HH$. The optimal constants $B_-\leq A_-<0<A_+\leq B_+$ satisfying Eq. \eqref{frm pos} are called the \emph{$J$-frame bounds of $\mc{F}$}. 
\end{defn}

In order to compute the $J$-frame bounds associated to a $J$-frame $\F=\{f_i\}_{i\in I}$,  consider the uniformly $J$-definite subspaces $\M_+$ and $\M_-$. Recall that $\F_+=\{f_i\}_{i\in I_+}$ is a frame for the Hilbert space $(\M_+,\K{\,}{\,})$. Then, if $G_+=G_{\M_+}|_{\M_+}\in GL(\M_+)$, the frame bounds for $\F_+$ are given by $A_+=\|(S_{G_+})^{-1}\|_+^{-1}$ and $B_+=\|S_{G_+}\|_+$, where $S_{G_+}=T_+T_+^*G_+$ is the frame operator of $\F_+$ and $\| f\|_+=\K{f}{f}^{1/2}=\|G_+^{1/2}f\|$, $f\in\M_+$, is the operator norm associated to the inner product $\K{\,}{\,}$. Therefore,
 \[
A_+=\|(S_{G_+})^{-1}\|_+^{-1}=\|G_+^{1/2}(T_+T_+^*G_+)^{-1}\|^{-1}= \|G_+^{-1/2}(T_+T_+^*)^{-1}\|^{-1},
 \]
and $B_+ = \|S_{G_+}\|_+= \|G_+^{1/2}T_+T_+^*G_+\|$. Analogously, it follows that $\F_-=\{f_i\}_{i\in I_-}$ is a frame for the Hilbert space $(\M_-, -\K{\,}{\,})$. So, the frame bounds for $\F_-$ are given by
\[
A_-= \|G_-^{-1/2}(T_-T_-^*)^{-1}\|^{-1} \ \ \ \text{and} \ \ \ B_-= \|G_-^{1/2}T_-T_-^*G_-\|,
\]
where $G_-=G_{\M_-}|_{\M_-}\in GL(\M_-)$. Thus, the $J$-frame bound associated to $\F$ can be fully characterized in terms of $T_\pm$ and the Gramian operators $G_{\M_\pm}$.

\subsection{Characterizing $J$-frames in terms of frame inequalities}

Given a Bessel family $\mc{F}=\{f_i\}_{i\in I}$ in a Krein space $\HH$, the inequalities:
\begin{equation}\label{frame ineq}
	A \ \K{f}{f} \leq \sum_{i\in I} |\K{f}{f_i}|^2 \leq B \ \K{f}{f} \ \ \ \text{for every $f\in \M=\ol{\Span\{f_i : i\in I\}}$},
\end{equation}
with $B\geq A>0$, ensure that $\M$ is a $J$-nonnegative subspace of $\HH$. However, they do not imply that $\M$ is uniformly $J$-positive, i.e. $(\M,\K{\,}{\,})$ is not necessarily a inner product space. See the example below. 

\begin{exmp}\label{ejem 2}

Consider again the Krein space $(\CC^3, \K{\,}{\,})$ as in Example \ref{ejem 1}. As it was mentioned before, $\M=\Span \{f_1=(1,0,1/\sqrt{2}),\,f_2=(0,1,1/\sqrt{2})\}$ is a $J$-nonnegative but not uniformly $J$-positive subspace of $\CC^3$. 

In this case, the orthogonal basis 
\[
v_1=(\tfrac{1}{2},\tfrac{1}{2},\tfrac{1}{\sqrt{2}})\, ,\ v_2=(\tfrac{1}{\sqrt{2}}\, , \ \tfrac{- 1}{\sqrt{2}}, 0)\ \text{ and }\ v_3=(\tfrac{1}{\sqrt{2}},\tfrac{1}{\sqrt{2}},-1),
\] 
is a basis of eigenvectors of $G_\M$, corresponding to the eigenvalues $\lambda_1=0$, $\lambda_2=1$ and $\lambda_3=0$, respectively. Moreover, $\M=\Span \{v_1,v_2\}$. 
Thus, if $f\in \mc M$ there exists $\alpha,\beta \in \CC$ such that $f=\alpha v_1+ \beta v_2$ and then, since $G_\M v_1=0\in \CC^3$, it is easy to see that 
\[
|\K{f}{f_1}|^2+|\K{f}{f_2}|^2=|\beta|^2 (|\PI{v_2}{f_1}|^2+|\PI{v_2}{f_2}|^2)=|\beta|^2=\K{f}{f}.
\] 
Therefore, Eq. \eqref{frame ineq} holds with $A=B=1$, but $\{f_1,\,f_2\}$ cannot be extended to a $J$-frame, since $\M$ is not a uniformly $J$-positive subspace.  
\end{exmp}

The next result gives a complete characterization of the families satisfying Eq. \eqref{frame ineq} for $B\geq A>0$. It is straightforward to formulate and prove analogues of all these assertions for a family satisfying Eq. \eqref{frame ineq} for negative constants $B\leq A<0$.

\begin{prop}\label{de desigualdades jframes}
Given a Bessel family  $\mc{F}=\{f_i\}_{i\in I}$ in a Krein space $\HH$, let $\mc{M}=\overline {\text{span}\{f_i\,:\ i\in I\}}$ and $\mc N=\M\cap\M^{\ort}$. If there exist constants $0<A\leq B$ such that
	\begin{equation}\label{lem frm pos}
	A \ \K{f}{f} \leq \sum_{i\in I} |\K{f}{f_i}|^2 \leq B \ \K{f}{f} \ \ \ \text{for every $f\in \M$},
	\end{equation}
then $\M\ominus \mc N$ is a (closed) uniformly $J$-positive subspace of $\mc M$. Moreover, if $\mc{F}$ is a frame for the Hilbert space $(\M, \PI{\,}{\,})$, the converse holds.
\end{prop}

\begin{pf*}{Proof}{}
First, suppose that there exist $0<A\leq B$ such that Eq. \eqref{lem frm pos} holds. So, $\M$ is a $J$-nonnegative subspace of $\HH$, or equivalently, $(\M,\K{\,}{\,})$ is a semi-inner product space. 
 
If $T\in L(\ell_2(I), \HH)$ is the synthesis operator of the Bessel sequence $\mc{F}$ and $C=\|T^*\|^2>0$, then $TT^* \leq C P_\M$. So, using Eq. \eqref{lem frm pos} it is easy to see that:
\begin{equation}\label{desigualdad PJP} 
A\PI{G_\M f}{f}  \leq \| T^\# (P_\M f)\|^2 = \PI{(P_\M JTT^*JP_\M ) f}{f} \leq C \PI{(G_\M)^2 f}{f}, \ \ f\in \HH.
\end{equation}
Thus, $0\leq G_\M\leq \frac{C}{A}\, (G_\M)^2$. Applying Douglas' theorem \cite{Douglas} it is easy to see that
\[
R((G_\M)^{1/2})\subseteq R(G_\M)\subseteq R((G_\M)^{1/2}).
\]
Moreover, it follows that $R(G_\M)$ is closed because $R(G_\M)=R((G_\M)^{1/2})$. 

Let $\M'=\M\ominus \N$ and notice that $\M'$ is a closed uniformly $J$-positive subspace of $\HH$. In fact, since $R(G_\M)$ is closed, there exists $\alpha >0$ such that
\[
\K{f}{f}=\PI{G_\M f}{f}=\|(G_\M)^{1/2} f\|^2 \geq \alpha \|f\|^2 \ \ \ \text{for every $f\in N(G_\M)^\bot=\M\ominus \N$}. 
\]
Conversely, suppose that $\mc F$ is a frame for $(\M,\PI{\,}{\,})$, i.e. there exist constants $B'\geq A'>0$ such that
\[
A' P_\M \leq TT^* \leq B' P_\M,
\]
where $T\in L(\ell_2(I),\M)$ is the synthesis operator of $\mc{F}$. If $\M'=\M\ominus \mc N$ is a uniformly $J$-positive subspace of $\HH$, then there exists $\alpha>0$ such that $\alpha P_{\M'} \leq G_{\M'} \leq P_{\M'}$. As a consequence of Douglas' theorem, $R((G_{\M'})^{1/2})=\M'=R(G_{\M'})$. Since $G_{\M}=G_{\M'}$ it is easy to see that
\[
A' (G_{\M})^2= A'(G_{\M'})^2 \leq P_\M JTT^*J P_\M \leq B' (G_{\M'})^2=B' (G_{\M})^2.
\]
Therefore, $R(P_\M JT)=R(G_{\M'})=R((G_{\M'})^{1/2})$, or equivalently, there exist $B\geq A>0$ such that
\[
A G_\M= A G_{\M'} \leq P_\M JTT^*J P_\M \leq B G_{\M'} = B G_{\M},
\]
i.e. $A \ \K{f}{f} \leq \sum_{i\in I} |\K{f}{f_i}|^2 \leq B \ \K{f}{f}$ for every $f\in \M$.
  \end{pf*}

\begin{thm}\label{thm J frame bounds}
Let $\mc{F}=\{f_i\}_{i\in I}$ be a frame for $\HH$. If $I_\pm=\{i\in I: \ \pm\K{f_i}{f_i}\geq 0\}$ and $\M_\pm=\ol{\Span\{f_i : i\in I_\pm\}}$ then, $\mc{F}$ is a $J$-frame if and only if $\M_\pm\cap \M_\pm^{\ort}=\{0\}$ and there exist constants $B_-\leq A_-<0<A_+ \leq B_+$ such that
\begin{equation}\label{eq J frame bounds}
	A_\pm \ \K{f}{f} \leq \sum_{i\in I_\pm} |\K{f}{f_i}|^2 \leq B_\pm \ \K{f}{f} \ \ \ \text{for every $f\in \M_\pm$}.
\end{equation}
\end{thm}

\begin{pf*}{Proof}{}
If $\mc{F}$ is a $J$-frame, the conditions on $\M_\pm$ follow by its definition and by Proposition \ref{jframes -> frames ineq}. Conversely, if $\M_+$ is $J$-non degenerated and there exist constants $0<A_+\leq B_+$ such that
\[
	A_+ \ \K{f}{f} \leq \sum_{i\in I_\pm} |\K{f}{f_i}|^2 \leq B_+ \ \K{f}{f} \ \ \ \text{for every $f\in \M_+$},
\]
then, by Proposition \ref{de desigualdades jframes}, $\M_+$ is a uniformly $J$-positive subspace of $\HH$. Therefore, there exist constants $0<A\leq B$ such that
\[
	A \ \|P_{\M_+} f\|^2 \leq \|T_+^\#P_{\M_+} f\|^2 \leq B \|P_{\M_+} f\|^2 \ \ \ \text{for every $f\in \HH$}.
\]
But these inequalities can be rewritten as 
\[
A\  P_{\M_+} \leq P_{\M_+} JT_+T_+^*JP_{\M_+} \leq B \ P_{\M_+}.
\]
Then, by Douglas' theorem, $R(P_{\M_+} JT_+)=R(P_{\M_+})=\M_+$. Furthermore, $P_{J(\M_+)}(R(T_+))=J(\M_+)$ because
\[
J(\M_+)=J(R(P_{\M_+} JT_+))=R((JP_{\M_+} J)T_+)=R(P_{J(\M_+)}T_+)=P_{J(\M_+)}(R(T_+)).
\]
Therefore, taking the counterimage of $P_{J(\M_+)}(R(T_+))$ by $P_{J(\M_+)}$, it follows that
\[
\HH=R(T_+) \dotplus J(\M_+)^\bot \subseteq \M_+ \dotplus \M_+^{\ort} =\HH.
\]
Thus, $R(T_+)=\M_+$ and $\mc{F}_+$ is a frame for $\M_+$. Analogously, $\mc{F}_-=\{f_i\}_{i\in I_-}$ is a frame for $\M_-$. Finally, since $\mc{F}$ is a frame for $\HH$,
\[
\HH= R(T)=R(T_+)+ R(T_-),
\]
which proves the maximality of $R(T_\pm)$. Thus, $\mc{F}$ is a $J$-frame for $\HH$.
  \end{pf*}

\subsection{A geometrical characterization of $J$-frames}

Let $\mc{F}=\{f_i\}_{i\in I}$ be a $J$-frame for $\HH$ and consider $\mc{F}=\mc{F}_+\cup \mc{F}_+$ the partition of $\mc{F}$ into $J$-positive and $J$-negative vectors. Moreover, let $\M_\pm$ be the (maximal) uniformly $J$-definite subspace of $\HH$ generated by $\mc{F}_\pm$. 

The aim of this section is to show that it is possible to bound the correlation between vectors in $\mc{F}_+$ (resp. $\mc{F}_-$) and vectors in the cone of neutral vectors $\mc{C}=\{ n\in \HH \, :\  \K{n}{n}=0\}$,  in a strong sense:
\begin{equation}\label{correlac1}
|\PI{f}{n}|\leq c_{\pm}\  \|f\|\ \|n\| \ , \ \ f\in \M_\pm \ , \ \ n\in \mc C\ ,
\end{equation} for some constants $\frac{\sqrt 2}{2}\leq c_\pm<1$. In order to make these ideas precise, consider the notion of minimal angle between a subspace $\M$ and the cone $\mc{C}$.

\begin{defn}\label{defi dist} 
Given a closed subspace  $\M$ of the Krein space $\HH$, consider 
\begin{equation}\label{eq distancia a neutros}
c_0(\M,\mc{C})=\sup\,\{|\PI{m}{n}|\,:\ m\in \M,\,n\in \mc{C}, \ \|n\|=\|m\|=1\}\ ,
\end{equation}
Then, there exists a unique $\theta(\M,\mc{C})\in [0,\frac{\pi}{4}]$ such that $\cos(\theta(\M,\mc{C}))=c_0(\M,\mc{C})$. In this case, $\theta(\M,\mc{C})$ is the \emph{minimal angle} between $\M$ and $\mc{C}$.
\end{defn}
\medskip

Observe that if the subspace $\M$ contains a non trivial $J$-neutral vector (e.g. if 	$\M$ is $J$-indefinite or $J$-semidefinite) then $c_0(\M,\mc{C})=1$, or equivalently, $\theta(\M,\mc{C})=0$. On the other hand, it will be shown that the minimal angle between a uniformly $J$-positive (resp. uniformly $J$-negative) subspace $\M$ and $\mc C$ is always bounded away from 0.

\begin{prop}\label{una cota sobre c0}
Let $\M$ be a $J$-semidefinite subspace of $\HH$ with definiteness bound $\alpha$. Then, 
\begin{equation}\label{eq angle}
c_0(\M,\mc{C})= \frac{1}{\sqrt{2}}\left( \sqrt{\frac{1+\alpha}{2}}+ \sqrt{\frac{1-\alpha}{2}} \,\right).
\end{equation}
In particular, $\M$ is uniformly $J$-definite if and only if $c_0(\M,\mc{C})<1$. 
\end{prop}

\begin{pf*}{Proof}{}
Let $\HH=\HH_+\oplus\HH_-$ be a fundamental decomposition of $\HH$ and suppose that $\M$ is a $J$-nonnegative subspace of $\HH$. 

Let $m\in \M$ with $\|m\|=1$. Then, there exist (unique) $m^\pm\in \HH_\pm$ such that $m=m^+ + m^-$. In this case, 
\begin{equation}\label{restrictions} 
1=\|m\|^2=\|m^+\|^2+\|m^-\|^2 \ \ \ \text{ and }\ \ \ \alpha\leq \K{m}{m}= \|m^+\|^2-\|m^-\|^2. 
\end{equation}

\noindent{\bf Claim:} Fixed $m\in \M$ with $\|m\|=1$, $\sup\,\{ |\PI{m}{n}|\, : \ n\in \mc{C}, \ \|n\|=1 \}=\frac{1}{\sqrt{2}}(\|m^+\|+\|m^-\|)$.

\medskip

Indeed, consider $n\in\mc{C}$ with $\|n\|=1$. Then, there exist (unique) $n_\pm\in \HH_\pm$ such that $n=n^+ + n^-$. In this case,
\[
0=\K{n}{n}=\|n^+\|^2-\|n^-\|^2 \ \ \ \text{ and }\ \ \ 1=\|n\|^2=\|n^+\|^2+\|n^-\|^2,
\] 
which imply that $\|n^+\|=\|n^-\|=\tfrac{1}{\sqrt{2}}$. Therefore, 
\[
|\PI{m}{n}|\leq |\PI{m^+}{n^+}|+ |\PI{m^-}{n^-}|\leq \frac{1}{\sqrt 2}\,(\|m^+\|+\|m^-\|).
\]
On the other hand, if $m^-\neq 0$ then let $n_m:=\frac{1}{\sqrt 2}(\frac{m^+}{\|m^+\|}+ \frac{m^-}{\|m^-\|} )$, otherwise consider $n_m=\tfrac{1}{\sqrt{2}}(m + z)$, with $z\in \HH_-$, $\|z\|=1$. Now, it is easy to see that $n_m\in \mc{C}$ and that $|\PI{m}{n_m}|= \frac{1}{\sqrt{2}}\,(\|m^+\|+\|m^-\|)$ which together with the previous facts prove the claim.

\medskip

Now, let $\M_1=\{ m=m^+ + m^-\in \M : \ m^\pm \in\HH_\pm, \ \|m\|=1\}$. Using the claim above it follows that
\begin{equation}\label{otra defi de la dist}
c_0(\M,\mc{C})=\frac{1}{\sqrt{2}} \ \sup_{m\in \M_1}  (\|m^+\|+\|m^-\|).
\end{equation}

If $\alpha=1$ then $\M$ is a subspace of $\HH_+$. Also, it is easy to see that $c_0(\M,\mc{C})=\tfrac{1}{\sqrt{2}}$. Thus, in this particular case, $c_0(\M,\mc{C})= \frac{1}{\sqrt{2}}\left( \sqrt{\frac{1+\alpha}{2}}+ \sqrt{\frac{1-\alpha}{2}} \,\right)$.

On the other hand, if $\alpha<1$, let $k_0\in\NN$ be such that $\frac{1-\alpha}{2}>\frac{1}{2k_0}$. Observe that, by the definition of the definiteness bound, for every integer $k\geq k_0$ there exists $m_k=m_k^+ + m_k^-\in\mc{M}_1$ such that $\alpha\leq \|m_k^+\|^2-\|m_k^-\|^2<\alpha +\frac{1}{k}$. Then, it follows that
\[
\alpha +1\leq 2\|m_k^+\|^2<\alpha +1+\frac{1}{k},
\]
or equivalently, $\sqrt{\frac{\alpha +1}{2}} \leq \|m_k^+\| < \sqrt{\frac{\alpha +1}{2}+\frac{1}{2k}}$. Moreover, $\|m_k^-\|=\sqrt{1-\|m_k^+\|^2}$ implies that
\[
\sqrt{\frac{1-\alpha}{2} -\frac{1}{2k}} < \|m_k^-\| \leq \sqrt{\frac{1-\alpha}{2}}.
\]
Therefore, for every integer $k\geq k_0$ there exists $m_k\in\mc{M}_1$ such that
\[
\sqrt{\frac{1-\alpha}{2} -\frac{1}{2k}}+\sqrt{\frac{\alpha +1}{2}} < \|m_k^+\|+\|m_k^-\| < \sqrt{\frac{\alpha +1}{2}+\frac{1}{2k}}+\sqrt{\frac{1-\alpha}{2}}.
\]
Thus, $c_0(\M,\mc{C})= \frac{1}{\sqrt{2}}\left( \sqrt{\frac{1+\alpha}{2}}+ \sqrt{\frac{1-\alpha}{2}} \,\right)$. 

Assume now that $\M$ is a $J$-nonpositive subspace of $(\HH, \K{\,}{\,})$ with definiteness bound $\alpha$, for $0\leq \alpha\leq 1$. Then, $\M$ is a $J$-nonnegative subspace of the antispace $(\HH,-\K{}{})$, with the same definiteness bound $\alpha$. Furthermore, the cone of $J$-neutral vectors for the antispace is the same as for the initial Krein space $(\HH, \K{\,}{\,})$. Therefore, we can apply the previous arguments and conclude that Eq.\eqref{eq angle} also holds for $J$-nonpositive subspaces.

Finally, the last assertion in the statement follows from the formula in Eq. \eqref{eq angle}.
  \end{pf*}

Let $\mc F$ be a $J$-frame for $\HH$ as above. Notice that the Eq. \eqref{correlac1} holds for some constant $\tfrac{\sqrt 2}{2}\leq c_\pm<1$ if and only if $c_0(\M_\pm,\mc{C})<1$, i.e. that the minimal angles $\theta(\M_\pm,\mc C)$ are bounded away from 0. This is intimately related with the fact that the aperture between the subspaces $\M_+$ (resp. $\M_-$) and $\HH_+$ (resp. $\HH_-$) is bounded away from $\frac{\pi}{4}$, whenever $\HH=\HH_+\oplus\HH_-$ is a fundamental decomposition.

\begin{rem}\label{rem apertur}
Given a Krein space $\HH$, fix a fundamental decomposition $\HH=\HH_+\oplus\HH_-$. Then, if $\M$ is a $J$-nonnegative subspace of $\HH$ the minimal angle between $\M$ and $\mc{C}$ is related with the \emph{aperture} $\Phi(\M,\HH_+)$ between the subspaces $\M$ and $\HH_+$, see \cite{Glazman} and Exercises 3--6 to \cite[Ch. 1, \S 8]{Iokhvidov}. In fact, if $K\in L(\HH_+,\HH_-)$ is the angular operator associated to $\M$ then, by \cite[Ch. 1, \S 8 Exercise 4]{Iokhvidov},
\[
\Phi(\M,\HH_+)= \frac{\|K\|}{\sqrt{1+ \|K\|^2}}.
\]
Also, if $\alpha$ is the definiteness bound of $\M$ then $\|K\|=\sqrt{\frac{1-\alpha}{1+\alpha}}$, see \cite[Ch. 1, Lemma 8.4]{Iokhvidov}. Therefore, $\Phi(\M,\HH_+)= \frac{\|K\|}{\sqrt{1+ \|K\|^2}}=\sqrt{\frac{1-\alpha}{2}}$.
Since $\Phi(\M,\HH_+)=\sin \vp(\M,\HH_+)$ for an angle $\vp (\M,\HH_+)\in [0,\tfrac{\pi}{4}]$ between $\M$ and $\HH_+$, it is easy to see that
\[
\cos \vp (\M,\HH_+)= \sqrt{1- \sin^2 \vp(\M,\HH_+)}=\sqrt{\frac{1+\alpha}{2}}.
\]
Therefore, if $\vp=\vp (\M,\HH_+)$,
\begin{eqnarray*}
\cos (\tfrac{\pi}{4}-\vp)= \frac{\sqrt{2}}{2}\left(\cos \vp + \sin\vp\right)= \frac{1}{\sqrt{2}}\left( \sqrt{\frac{1+\alpha}{2}}+ \sqrt{\frac{1-\alpha}{2}} \,\right)=c_0(\M,\mc{C}),
\end{eqnarray*}
i.e. $\vp (\M,\HH_+) + \theta(\M,\mc{C})=\frac{\pi}{4}$.  
\end{rem}

The following result shows that, given a frame $\mc{F}=\{f_i\}_{i\in I}$ for $\HH$, the positivity of the angles $\theta(\M_\pm,\mc C)$ characterize it as a $J$-frame for $\HH$.

\begin{prop}
Let $\F=\{f_i\}_{i\in I}$ be a frame for a Krein space $\HH$. Then, $\F$ is a $J$-frame for $\HH$ if and only if there exists a partition $I=I_1\cup I_2$ such that
\begin{equation}\label{ecua proteo} 
\theta(\M_j,\mc C)>0 \ \ \ \text{for $j=1,2$},
\end{equation}
where $\M_j=\ol{\Span\{f_i : i\in I_j\}}$. 
\end{prop}

\begin{pf*}{Proof}{}
If we assume that $\F$ is a $J$-frame then, consider $I_\pm$ and $\M_\pm$ as usual. Then $I=I_+\cup I_-$ is a partition of $I$ into disjoint sets and $\M_\pm$ are uniformly $J$-definite subspaces associated to $\F$. Hence, by Proposition \ref{una cota sobre c0}, we see that Eq. \eqref{ecua proteo} holds in this case.

Conversely, assume that there exists a partition of $I$ with the properties above. Notice that Proposition \ref{una cota sobre c0} implies that $\M_j$ is a uniformly $J$-definite subspace of $\HH$, for $j=1,2$. On the other hand, since $\F$ is a frame, $\HH\subseteq \M_1 + \M_2$. 
Therefore, $\M_1$ and $\M_2$ have different positivity and they are maximal uniformly $J$-definite subspaces. Suppose that $\M_1$ is uniformly $J$-positive and $\M_2$ is uniformly $J$-negative. 

Then, consider the orthogonal projection $P_j\in L(\ell_2(I))$ onto the subspace $\ell_2(I_j)$, for $j=1,2$. If $\F=\{f_i\}_{i\in I}$ is a frame for $\HH$, its synthesis operator $T\in L(\ell_2(I),\HH)$ is surjective. Therefore,
\[
R(T_1)\dotplus R(T_2) = R(T) =\HH,
\]
where $T_j=TP_j$, for $j=1,2$. Then, it is easy to see that $R(T_j)=\M_j$ for $j=1,2$ and $\F$ is a $J$-frame for $\HH$.
  \end{pf*}

\begin{rem}
Let $(\HH, \PI{\,}{\,})$ be a separable Hilbert space that models a signal space in which is considered a linear (robust and stable) encoding-decoding scheme for certain measurements, i.e. consider a (redundant) frame $\mc{G}=\{g_i\}_{i\in K}$ for $\HH$. 

Assume that the measurements of $x\in\HH$ are given by $y_1=P x$ and $y_2=(I-P) x$, where $P \in L(\HH)$ is an orthogonal projection (for instance, $P$ and $I-P$ are low pass and high pass filters, respectively). Suppose that the signals having the same energy in $R(P)$ and $R(I-P)=N(P)$ (i.e. signals $x\in\HH$ such that $\|y_1\|^2=\|y_2\|^2$) are considered disturbances, see e.g. \cite{Buades, Masoum}. 

Notice that, sampling the measurements $y_1, y_2$ with the frame $\mc{G}$ is the same as sampling $y=(y_1,y_2)\in\HH\x\HH$ with the frame $\F=\{f_i\}_{i\in I}=\{(g_i,0)\}_{i\in K}\cup\{(0,g_i)\}_{i\in K}$ for $\HH\x\HH$.

It is easy to see that, the space $\KK=\HH\x\HH$ with the indefinite product $\K{y}{z}=\PI{y_1}{z_1}-\PI{y_2}{z_2}$ is a Krein space, where $y=(y_1,y_2), z=(z_1,z_2)\in\KK$ are the measurements of signals in $\HH$. Observe that the set of disturbances is characterized as the set of $J$-neutral vectors $\mc{C}$ of $\KK$.  

Also, notice that $\F$ is a $J$-frame for $\KK$. Hence, the (sampling) vectors of the frame $\F$ are away from the disturbances set $\mc{C}$. 

Now, consider any (redundant) $J$-frame $\mc{F}=\{f_i\}_{i\in I}$ for $(\KK,\K{\,}{\,})$. As usual, denote $\M_+$ and $\M_-$ the maximal uniformly $J$-definite subspaces generated by $\F$. Since $\M_\pm$ is uniformly $J$-definite, Proposition \ref{una cota sobre c0} shows that $c_0(\M_\pm,\mc{C})<1$, which is a bound for the correlation between the sampling vectors in $\F$ and the distrubances of $\mc{C}$ because
\begin{equation}\label{desi conoc1}
|\PI{f_i}{n}|\leq c_0(\M_\pm,\mc{C}) \, \|f_i\|\, \|n\|\ \ \ \text{whenever $i\in I_\pm$ and $n\in \mc{C}$}.
\end{equation} 
That is, $J$-frames provide a class of frames for $\KK$ with the desired properties. Moreover, later in Proposition \ref{J dual}, it will be shown that the $J$-frame $\F$ admits a (canonical) dual $J$-frame that induces a linear (indefinite) stable and redundant encoding-decoding scheme in which the correlation between
both the sampling and reconstructing vectors and the cone of neutral vectors is bounded from above. These remarks provide a quantitative measure of the advantage of considering $J$-frames with respect to usual frames in this setting. \end{rem}

\section{On the synthesis operator of a $J$-frame}

If $\mc{F}$ is a $J$-frame with synthesis operator $T$, then $QT=T_+=TP_+$, where $Q=P_{\M_+//\M_-}$. Therefore, 
\[
Q=QTT^\dag=TP_+T^\dag.
\]
So, given a surjective operator $T:\ell_2(I)\ra \HH$, the idempotency of $TP_+T^\dag$ is a necessary condition for $T$ to be the synthesis operator of a $J$-frame.

\begin{lem}\label{NTbot}
Let $T\in L(\ell_2(I),\HH)$ be surjective. Suppose that $P_\St$ is the orthogonal projection onto a closed subspace $\St$ of $\ell_2(I)$ such that $c(\St, N(T)^\bot)<1$. Then, $TP_\St T^\dag$ is a projection if and only if 
\[
N(T)= \St\cap N(T) \oplus \St^\bot\cap N(T).
\]
\end{lem}

\begin{pf*}{Proof}{}
Suppose that $Q=TP_\St T^\dag$ is a projection. Then, if $P=P_{N(T)^\bot}$, $E=PP_\St P$ is an orthogonal projection because it is selfadjoint and
\[
E^2=(PP_\St P)^2=PP_\St P P_\St P=T^\dag(TP_\St T^\dag)^2T=T^\dag(TP_\St T^\dag)T=P P_\St P=E.
\]
Therefore, $(PP_\St)^k= E^{k-1}P_\St=EP_\St=(PP_\St)^2$ for every $k\geq 2$. So, by \cite[Lemma 18]{Deu}, 
\[
P P_\St=P_\St \wedge P=P_\St P.
\]
Then, since $P_\St$ and $P$ commute, it follows that $N(T)= \St\cap N(T) \oplus \St^\bot\cap N(T)$ (see \cite[Lemma 9]{Deu}).

Conversely, suppose that $N(T)= \St\cap N(T) \oplus \St^\bot\cap N(T)$. Then, $P_\St$ and $P$ commute and
\[
(TP_\St T^\dag)^2=TP_\St(T^\dag T)P_\St T^\dag=TP_\St P P_\St T^\dag =TP P_\St T^\dag =TP_\St T^\dag.  
\]
\end{pf*}

Hereafter consider the set of possible decompositions of $\HH$ as a (direct) sum of a pair of maximal uniformly definite subspaces, or equivalently, the associated set of projections:
\[
\Q =\{Q\in L(\HH) : Q^2=Q, \, \text{$R(Q)$ is uniformly $J$-positive and $N(Q)$ is uniformly $J$-negative} \}.
\]

\begin{prop}\label{N T}
Let $T\in L(\ell_2(I),\HH)$ be surjective. Then, $T$ is the synthesis operator of a $J$-frame if and only if there exists $I_+\subset I$ such that $\ell_2(I_+)$ (as a subspace of $\ell_2(I)$) satisfies $c(N(T)^\bot, \ell_2(I_+))<1$ and
\[
TP_+T^\dag \in \Q,
\]
where $P_+\in L(\ell_2(I))$ is the orthogonal projection onto $\ell_2(I_+)$.
\end{prop}

\begin{pf*}{Proof}{}
If $T$ is the synthesis operator of a $J$-frame, the existence of such a subset $I_+$ has already been discussed before.

Conversely, suppose that  there exists such a subset $I_+$ of $I$. Then, since $c(N(T)^\bot, \ell_2(I_+))<1$ and $Q=TP_+T^\dag\in \Q$, it follows from Lemma \ref{NTbot} that $P_+$ and $P=P_{N(T)^\bot}$ commute. Therefore, 
\[
QT=TP_+P=TP P_+= TP_+,
\]
and $(I-Q)T=T(I-P_+)$. Hence, $R(TP_+)=R(Q)$ is (maximal) uniformly $J$-positive and $R(T(I-P_+))=N(Q)$ is (maximal) uniformly $J$-negative. Therefore $\F=\{Te_i\}_{i\in I}$ is by definition a $J$-frame for $\HH$.
  
\end{pf*}

\begin{thm}
Given a surjective operator $T\in L(\ell_2(I), \HH)$, the following conditions are equivalent:
\begin{enumerate}
	\item There exists $U\in \U(\ell_2(I))$ such that $TU$ is the synthesis operator of a $J$-frame.
	\item There exists $Q\in \Q$ 	such that 
	\begin{equation}\label{QTT(I-Q)=0}
	QTT^*(I-Q)^*=0.
	\end{equation}
	\item There exist closed range operators $T_1, T_2\in L(\ell_2(I),\HH)$ such that $T=T_1 + T_2$, $R(T_1)$ is uniformly $J$-positive, $R(T_2)$ is uniformly $J$-negative and $T_1T_2^*=T_2T_1^*=0$.
\end{enumerate}
\end{thm}

\begin{pf*}{Proof}{}
\noi {\it 1.} $\Rightarrow$ {\it 2.}: \ Suppose that there exists $U\in \U(\ell_2(I))$ such that $V=TU$ is the synthesis operator of a $J$-frame. If $I_\pm=\{i\in I: \pm\K{Ve_i}{Ve_i} >0\}$ and $P_\pm\in L(\ell_2(I))$ is the orthogonal projection onto $\ell_2(I_\pm)$, define $V_\pm=VP_\pm$. Then, $V=V_+ + V_-$ and $\M_\pm=R(V_\pm)$ is a maximal uniformly $J$-definite subspace. So, considering $Q=P_{\M_+//\M_-}\in\Q$, it is easy to see that $QV=V_+$, $(I-Q)V=V_-$ and
\[
QTT^*(I-Q)^*=QVV^*(I-Q)^*=V_+V_-^*=VP_+P_-V^*=0.
\]
\noi {\it 2.} $\Rightarrow$ {\it 3.}: \ Suppose that there exists $Q\in \Q$ such that $QTT^*(I-Q)^*=0$. Defining $T_1=QT$ and $T_2=(I-Q)T$, it follows that $T=T_1 + T_2$, $R(T_1)=R(Q)$ is uniformly $J$-positive, $R(T_2)=N(Q)$ is uniformly $J$-negative and
\[
T_1T_2^*=T_2T_1^*=0,
\]
because Eq. \eqref{QTT(I-Q)=0} says that $R(T_2^*)=R(T^*(I-Q)^*)\subseteq N(QT)= N(T_1)$.

\noi {\it 3.} $\Rightarrow$ {\it 1.}: \ If there exist closed range operators $T_1, T_2\in L(\ell_2(I),\HH)$ satisfying the conditions of item {\it 3.}, notice that $T_1T_2^*=0$ implies that $N(T_2)^\bot\subseteq N(T_1)$, or equivalently, $N(T_1)^\bot\subseteq N(T_2)$. 

Consider the projection $Q=P_{R(T_1)//R(T_2)}\in\Q$ and notice that $QT=T_1$ and $(I-Q)T=T_2$. If $\mc{B}_1=\{u_i\}_{i\in I_1}$ is an orthonormal basis of $N(T_1)^\bot$, consider the family $\{f_i^+\}_{i\in I_1}$ in $\HH$ given by $f_i^+=Tu_i$. But, if $i\in I_1$, 
\[
f_i^+=QTu_i + (I-Q)Tu_i=T_1u_i\in R(T_1),
\]
because $u_i\in N(T_1)^\bot\subseteq N(T_2)$. Therefore, $\{f_i^+\}_{i\in I_1}\subseteq R(T_1)$. Since $T_1$ is an isomorphism between $N(T_1)^\bot$ and $R(T_1)$, it follows that $R(T_1)=\ol{\Span \{f_i^-\}_{i\in I_1}}$. 

Analogously, if $\mc{B}_2=\{b_i\}_{i\in I_2}$ is an orthonormal basis of $N(T_1)$, the family $\{f_i^-\}_{i\in I_2}$ defined by $f_i^-=Tb_i$ ($i\in I_2$) lies in $R(T_2)$.
Since $T_2$ is an isomorphism between $N(T_2)^\bot$ and $R(T_2)$, it follows that 
\[
R(T_2)= T_2(N(T_1))\subseteq \ol{\Span\{f_i^-\}_{i\in I_2}}\subseteq R(T_2).
\] 
Finally, consider $U\in \U(\ell_2(I))$ which turns the standard orthonormal basis $\{e_i\}_{i\in I}$ into $\mc{B}_1\cup\mc{B}_2$. Then, if $V=TU$ and $\mc{F}=\{Ve_i\}_{i\in I}=\{f_i^+\}_{i\in I_1}\cup\{f_i^-\}_{i\in I_2}$, it is easy to see that 
\[
I_+=\{i\in I : \K{Ve_i}{Ve_i}>0\}=I_1 \ \ \ \text{and} \ \ \ I_-=\{i\in I : \K{Ve_i}{Ve_i}<0\}=I_2.
\]
So, $R(V_+)=R(T_1)$ is maximal uniformly $J$-positive and $R(V_-)=R(T_2)$ is maximal uniformly $J$-negative. Therefore, $\mc{F}$ is a $J$-frame for $\HH$ with synthesis operator $V=TU$.
  \end{pf*}

\section{The $J$-frame operator}

\begin{defn}  
Given a $J$-frame $\mc{F}=\{f_i\}_{i\in I}$, the \emph{$J$-frame operator} $S:\HH\ra \HH$ is defined by
\[
Sf=\sum_{i\in I}\s_i \K{f}{f_i} f_i, \ \ \ \text{for every $f\in \HH$},
\]
where $\s_i=\sgn(\K{f_i}{f_i})$.
\end{defn}

\noindent The following proposition compiles some basic properties of the $J$-frame operator. 

\begin{prop}\label{props S}
Let $\mc{F}=\{f_i\}_{i\in I}$ be a $J$-frame with synthesis operator $T\in L(\ell_2(I), \HH)$. Then, its $J$-frame operator $S\in L(\HH)$ satisfies:
\begin{enumerate}
	\item $S=TT^\#$; 
	\item $S=S_+ - S_-$, where $S_+:=T_+T_+^\#$ and $S_-:=-T_-T_-^\#$ are $J$-positive operators;
	\item $S$ is an invertible $J$-selfadjoint operator;
	\item $\ind_\pm(S)=\dim \HH_\pm$, where $\ind_\pm(S)$ are the indices of $S$.
\end{enumerate}
\end{prop}

\begin{pf*}{Proof}{}
If $\mc{F}=\{f_i\}_{i\in I}$ is a $J$-frame with synthesis operator $T\in L(\ell_2(I), \HH)$, then $T^\#f=\sum_{i\in I} \s_i \K{f}{f_i}e_i$ for $f\in \HH$. So,
\[
TT^\#f= T\left(\sum_{i\in I} \s_i \K{f}{f_i}e_i\right)= \sum_{i\in I}\s_i \K{f}{f_i} f_i=Sf, \ \ \ \text{for every $f\in \HH$}.
\]
Furthermore, if $I_\pm=\{i\in I: \pm\K{f_i}{f_i}>0\}$, consider $T_\pm=TP_\pm$ as usual. Then,
	\[
	TT^\# = (T_+ + T_-)(T_+ + T_-)^\#= T_+T_+^\# + T_-T_-^\#= T_+T_+^\# - (-T_-T_-^\#),
	\]
because $T_+T_-^\#=T_-T_+^\#=0$. Therefore, $S=S_+ - S_-$ if $S_\pm:= \pm T_\pm T_\pm^\#$.
Notice that $S_\pm$ is a $J$-positive operator because
	\[
	S_\pm= \pm T_\pm T_\pm^\#= \pm T_\pm J_2T_\pm^* J=  T_\pm T_\pm^*J. 
	\]

To prove the invertibility of $S$ observe that, if $Sf=0$ then $S_+f=S_-f$. But $R(S_+)\cap R(S_-)\subseteq R(T_+)\cap R(T_-)=\{0\}$. Thus, $S$ is injective. On the other hand, $R(S)=S(\M_+^{\ort}) + S(\M_-^{\ort})$ because $\HH=\M_+^{\ort}\dotplus \M_-^{\ort}$. But it is easy to see that $\M_\pm^{\ort}\subseteq N(S_\pm)$. So, $S(\M_\pm^{\ort})=S_\mp(\M_\pm^{\ort})$ and $R(S)= S_-(\M_+^{\ort})+ S_+(\M_-^{\ort})= R(S_-) + R(S_+)= \M_+ + \M_-=\HH$. Therefore, $S$ is invertible.

Finally, the identities $\ind_\pm(S)=\dim \HH_\pm$ follow from the indices definition. Recall that if $A\in L(\HH)$ is a $J$-selfadjoint operator, $\ind_+(A)$ is the supremum of all positive integers $r$ such that there exists a positive invertible matrix of the form $(\K{Ax_j}{x_k})_{j,k=1,\ldots,r}$, where $x_1,\ldots, x_r\in\HH$ (if no such $r$ exists, $\ind_-(A) = 0$). Similarly, $\ind_-(A) = \ind_+(-A)$ is the supremum of all positive integers $m$ such that there exists a negative invertible matrix of the form $(\K{Ay_j}{y_k})_{j,k=1,\ldots,m}$, where $y_1,\ldots, y_m\in\HH$, see \cite{Dritschel 1}. 
  \end{pf*}
 
\begin{cor}\label{S+ y S-}
Let $\mc{F}=\{f_i\}_{i\in I}$ be a $J$-frame for $\HH$ with $J$-frame operator $S\in L(\HH)$. Then, $R(S_\pm)=\M_\pm$ and $N(S_\pm)=\M_\pm^{\ort}$. Furthermore, if $Q=P_{\M_+//\M_-}$,
\begin{equation}
S_+=QSQ^\# \ \ \ \text{and} \ \ \ S_-=-(I-Q)S(I-Q)^\#.
\end{equation}
\end{cor} 
 
\begin{pf*}{Proof}{}
Recall that $S_+:= T_+ T_+^\#=T_+(J_2 T_+^* J)= T_+T_+^* J$. Then, $R(S_+)=R(T_+T_+^*J)=R(T_+T_+^*)=R(T_+)=\M_+$ because $R(T_+)$ is closed. Since $S_+$ is $J$-selfadjoint, it follows that $N(S_+)=R(S_+)^{\ort}=\M_+^{\ort}$.
Analogously, $R(S_-)=\M_-$ and $N(S_-)=\M_-^{\ort}$.

Since $S=S_+ - S_-$, if $Q=P_{\M_+//\M_-}$ then
\[
QS=Q(S_+ - S_-)=S_+,
\]
by the characterization of the range and nullspace of $S_+$. Therefore, $SQ^\#=QS=QSQ^\#$.
Analogously, $S(I-Q)^\#=(I-Q)S=(I-Q)S(I-Q)^\#$.
  \end{pf*}
 
The above corollary states that $S$ is the diagonal block operator matrix
\begin{equation}\label{S en bloques}
S=\matriz{S_+}{0}{0}{-S_-},
\end{equation}
according to the (oblique) decompositions $\HH=\M_-^{\ort} \dotplus \M_+^{\ort}$ and $\HH=\M_+\dotplus \M_-$ of the domain and codomain of $S$, respectively.

\subsection{The indefinite reconstruction formula associated to a $J$-frame}

Given a $J$-frame $\mc{F}=\{f_i\}_{i\in I}$ with synthesis operator $T$, there is a duality between $\mc{F}$ and the frame $\mc{G}=\{g_i\}_{i\in I}$ given by $g_i=S^{-1}f_i$: if $f\in\HH$,
\begin{eqnarray*}
f= S S^{-1}f &=& TT^\#(S^{-1}f)=T\left(\sum_{i\in I} \s_i\K{S^{-1}f}{f_i}e_i \right)= \sum_{i\in I} \s_i\K{S^{-1}f}{f_i}f_i = \sum_{i\in I} \s_i\K{f}{S^{-1}f_i}f_i.
\end{eqnarray*}
Analogously, 
\[
f= S^{-1}S f = S^{-1}(TT^\#f) =S^{-1} \left(\sum_{i\in I} \s_i\K{f}{f_i} f_i\right) = \sum_{i\in I} \s_i\K{f}{f_i} S^{-1}f_i.
\]
Therefore, for every $f\in \HH$, there is an \emph{indefinite reconstruction formula} associated to $\F$:
\begin{equation}\label{Jreconst}
	f= \sum_{i\in I} \s_i\K{f}{g_i}f_i = \sum_{i\in I} \s_i\K{f}{f_i} g_i.
\end{equation}
The following question arises  naturally: is $\mc{G}=\{S^{-1}f_i\}_{i\in I}$ also a $J$-frame for $\HH$?

\begin{prop}\label{J dual}
If $\mc{F}=\{f_i\}_{i\in I}$ is a $J$-frame for a Krein space $\HH$ with $J$-frame operator $S$, then $\mc{G}=\{S^{-1}f_i\}_{i\in I}$ is also a $J$-frame for $\HH$.
\end{prop}

\begin{pf*}{Proof}{}
Given a $J$-frame $\mc{F}=\{f_i\}_{i\in I}$ for $\HH$ with $J$-frame operator $S$, observe that the synthesis operator of  $\mc{G}=\{S^{-1}f_i\}_{i\in I}$ is $V:=S^{-1}T\in L(\ell_2(I),\HH)$. Furthermore, by Corollary \ref{S+ y S-}, $S(\M_\mp^{\ort})=\M_\pm$. Then, $S^{-1}(\M_\pm)= \M_\mp^{\ort}$ and it follows that
\[
\K{S^{-1}f_i}{S^{-1}f_i}>0 \ \ \ \text{if and only if} \ \ \ \K{f_i}{f_i}>0. 
\]
Thus, $V_\pm=VP_\pm=S^{-1}T_\pm$ and $R(V_+)$ (resp. $R(V_-)$) is a maximal uniformly $J$-positive (resp. $J$-negative) subspace of $\HH$. So, $\mc{G}$ is a $J$-frame for $\HH$.
  \end{pf*}

If $\mc{F}=\{f_i\}_{i\in I}$ is a frame for a Hilbert space $\HH$ with synthesis operator $T \in L(\ell_2(I), \HH$), then the family $\{(TT^*)^{-1}f_i\}_{i\in I}$ is called the \emph{canonical dual frame} because it is a dual frame for $\F$ (see Eq. \eqref{dual}) and it has the following optimal property: Given $f\in \HH$,
\begin{equation}\label{canonico}
\sum_{i\in I} |\PI{f}{(TT^*)^{-1}f_i}|^2 \leq \sum_{i\in I} |c_i|^2, \ \ \ \text{whenever} \ \ \ f=\sum_{i\in I} c_i f_i,
\end{equation}
for a family $(c_i)_{i\in I}\in \ell_2(I)$. In other words, the above representation has the smallest $\ell_2$-norm among the admissible frame coefficients representing $f$ (see \cite{DS}).

\begin{rem}
If $\F=\{f_i\}_{i\in I}$ is a $J$-frame for $\HH$ then $\F_\pm=\{f_i\}_{i\in I_\pm}$ is a frame for the Hilbert space $(\M_\pm,\pm\K{\,}{\,})$. Furthermore, the frame operator associated to $\F_+$ is $S_+=T_+T_+^\#$ and its canonical dual frame is given by $\mc{G}_+=\{S_+^{-1}f_i\}_{i\in I_+}$. Analogously, the frame operator associated to $\F_-$ is $S_-=-T_-T_-^\#$ and its canonical dual frame is given by $\mc{G}_-=\{S_-^{-1}f_i\}_{i\in I_-}$.

Then, since $\HH=\M_+ \dotplus \M_-$, $\HH$ can be seen as the (outer) direct sum of the Hilbert spaces $(\M_+, \K{\,}{\,})$ and $(\M_-,-\K{\,}{\,})$, i.e. the inner product given by
\[
\PI{f}{g}_\F=\K{f_+}{g_+} - \K{f_-}{g_-}, \ \ \ f=f_+ + f_-, \ g=g_+ + g_-, \ \ \ f_+,g_+\in\M_+, \ f_-,g_-\in\M_-,
\]
turns $(\HH,\PI{\,}{\,}_\F)$ into a Hilbert space and the projection $Q=P_{\M_+//\M_-}$ is selfadjoint in this Hilbert space. So, if $f\in \HH$,
\begin{eqnarray*}
	\sum_{i\in I} |\K{f}{S^{-1}f_i}|^2 &=& \sum_{i\in I_+} |\K{Qf}{S_+^{-1}f_i}|^2 + \sum_{i\in I_-} |\K{(I-Q)f}{S_-^{-1}f_i}|^2 \leq \\ &\leq& \sum_{i\in I_+} |c_i^+|^2 + \sum_{i\in I_-} |c_i^-|^2,
\end{eqnarray*}
whenever $f_+=Qf= \sum_{i\in I_+} c_i^+ f_i$ and $f_-=(I-Q)f= \sum_{i\in I_-} c_i^- f_i$, for families $(c_i^\pm)_{i\in I_\pm}\in \ell_2(I_\pm)$. Therefore,
\[
\sum_{i\in I} |\K{f}{S^{-1}f_i}|^2 \leq \sum_{i\in I} |c_i|^2,
\]
whenever $f= \sum_{i\in I} c_i f_i$ for some $(c_i)_{i\in I}\in \ell_2(I)$. In other words, the $J$-frame $\mc{G}=\{S^{-1}f_i\}_{i\in I}$ is the canonical dual frame of $\F$ in the Hilbert space $(\HH, \PI{\,}{\,}_\F)$.
\end{rem}

\subsection{Characterizing the $J$-frame operators}

In a Hilbert space $\HH$, it is well known that every positive invertible operator  $S\in L(\HH)$ can be realized as the frame operator of a frame $\mc F=\{f_i\}_{i\in I}$ for $\HH$, see \cite{HanLarson}. Indeed, if $\mc B=\{x_i\}_{i\in I}$ is an orthonormal basis of $\HH$, consider $T:\ell_2(I)\rightarrow \HH$ given by $Te_i= S^{1/2} x_i$ for $i\in I$. Then, for every $f\in \HH$, 
\[
TT^*f=\sum_{i\in I} \PI{f}{S^{1/2} x_i}\, S^{1/2} x_i=S^{1/2} \left( \sum_{i\in I} \PI{S^{1/2} f}{ x_i} \,  x_i\right )= S f 
\] 
Therefore, $\mc F=\{S^{1/2} x_i\}_{i\in I}$ is a frame for $\HH$ and its frame operator is given by $S$.

\medskip

The following paragraphs are devoted to characterize the set of $J$-frame operators.

\begin{thm}\label{carac S}
Let $S\in GL(\HH)$ be a $J$-selfadjoint operator acting on a Krein space $\HH$ with fundamental symmetry $J$. Then, the following conditions are equivalent: 
\begin{enumerate}
	\item $S$ is a $J$-frame operator, i.e. there exists a $J$-frame $\mc{F}$ with synthesis operator $T$ such that $S=TT^\#$.
	\item There exists a projection $Q\in \Q$ such that $QS$ is $J$-positive and $(I-Q)S$ is $J$-negative. 
	\item There exist $J$-positive operators $S_1,S_2\in L(\HH)$ such that $S=S_1 - S_2$ and $R(S_1)$ (resp. $R(S_2)$) is a uniformly $J$-positive (resp. $J$-negative) subspace of $\HH$.
\end{enumerate}
\end{thm}

\begin{pf*}{Proof}{}
{\it 1.} $\ra$ {\it 2.} follows from Proposition \ref{props S} and Corollary \ref{S+ y S-}.

{\it 2.} $\ra$ {\it 3.}: If there exists a projection $Q\in \Q$ such that $QS$ is $J$-positive and $(I-Q)S$ is $J$-negative, consider the $J$-positive operators $S_1=QS$ and $S_2=-(I-Q)S$. Then, $S=S_1-S_2$ and, by hypothesis, $R(S_1)=R(Q)$ is uniformly $J$-positive and $R(S_2)=R(I-Q)=N(Q)$ is uniformly $J$-negative.

{\it 3.} $\ra$ {\it 1.}: Suppose that there exist $J$-positive operators $S_1,S_2\in L(\HH)$ such that $S=S_1 - S_2$ and $R(S_1)$ (resp. $R(S_2)$) is a uniformly $J$-positive (resp. $J$-negative) subspace of $\HH$. Denoting $\KK_j=R(S_j)$ for $j=1,2$, observe that $A_j=S_j J|_{\KK_j}\in GL(\KK_j)^+$. Therefore, there exists a frame $\mc{F}_j=\{f_i\}_{i\in I_j}\subset \KK_j$ for $\KK_j$ such that $A_j=T_jT_j^*$ if $T_j\in L( \ell_2(I_1), \KK_j)$ is the synthesis operator of $\mc{F}_j$, for $j=1,2$. 

Then, consider $\ell_2(I):=\ell_2(I_1)\oplus\ell_2(I_2)$ and $T\in L(\ell_2(I), \HH)$ given by
\[
Tx= T_1x_1 + T_2x_2, \ \ \ \ \ \text{if $x\in \ell_2(I)$, $x=x_1 + x_2$, $x_j\in \ell_2(I_j)$ for $j=1,2$.}
\]
It is easy to see that $T$ is the synthesis operator of the frame $\mc{F}=\mc{F}_1\cup\mc{F}_2$. Furthermore $\mc{F}$ is a $J$-frame such that $I_+=I_1$ and $I_-=I_2$. 

Finally, endow $\ell_2(I)$ with the indefinite inner product defined by the diagonal operator $J_2\in L(\ell_2(I))$ given by
\[
J_2\, e_i= \s_i\, e_i, 
\]
where $\s_i=1$ if $i\in I_1$ and $\s_i=-1$ if $i\in I_2$. Notice that $T_1J_2=T_1$ and $T_2J_2=-T_2$. Furthermore, $T_1T_2^*=T_2T_1^*=0$ because $R(T_2^*)=N(T_2)^\bot\subseteq \ell_2(I_1)^\bot=\ell_2(I_2)\subseteq N(T_1)$. Thus,
\[
TT^\#=TJ_2T^*J=(T_1 + T_2)(T_1^*-T_2^*)J=T_1T_1^*J-T_2T_2^*J=A_1J-A_2J=S_1-S_2=S.   
\]
\end{pf*}

Given a $J$-frame  $\mc{F}=\{f_i\}_{i\in I}$ for $\HH$ with $J$-frame operator $S\in L(\HH)$, it follows from Corollary \ref{S+ y S-} that
\begin{equation}
S(\M_-^{\ort})=\M_+ \ \ \ \text{and} \ \ \ S(\M_+^{\ort})=\M_-.
\end{equation}
i.e. $S$ maps a maximal uniformly $J$-positive (resp. $J$-negative) subspace into another maximal uniformly $J$-positive (resp. $J$-negative) subspace. The next proposition shows under which hypotheses the converse holds.	

\begin{prop}\label{manda pos en pos}
Let $S\in GL(\HH)$ be a $J$-selfadjoint operator. Then, $S$ is a $J$-frame operator if and only if the following conditions hold:
\begin{enumerate}
	\item there exists a maximal uniformly $J$-positive subspace $\T$ of $\HH$ such that $S(\T)$ is also maximal uniformly $J$-positive;
	\item $\K{Sf}{f}\geq 0$ for every $f\in\T$;
	\item $\K{Sg}{g}\leq 0$ for every $g\in S(\T)^{\ort}$.
\end{enumerate}
\end{prop}

\begin{pf*}{Proof}{}
If $S$ is a $J$-frame operator, consider $\T=\M_-^{\ort}$ which is a maximal uniformly $J$-positive subspace $\T$ of $\HH$. Then, $S(\T)=\M_+$ is also maximal uniformly $J$-positive. Furthermore,
\[
\K{Sf}{f}=\K{SQ^\#f}{Q^\#f}=\K{QSQ^\#f}{f}=\K{S_+f}{f}\geq 0 \ \ \ \text{for every $f\in\T$},
\]
where $Q=P_{\M_+//\M_-}$. Also, $S(\T)^{\ort}=\M_+^{\ort}=N(Q^\#)=R((I-Q)^\#)$. So, 
\[
\K{Sg}{g}=\K{S(I-Q)^\#g}{(I-Q)^\#g}=\K{(I-Q)S(I-Q)^\#g}{g}=\K{-S_-g}{g}\leq 0 \ \ \ \text{for every $g\in S(\T)^{\ort}$}.
\]

Conversely, suppose that there exists a maximal uniformly $J$-positive subspace $\T$ satisfying the hypothesis. Let $\M=S(\T)$, which is maximal uniformly $J$-positive. Then, consider $Q=P_{\M//\T^{\ort}}$. It is well defined because $\T^{\ort}$ is maximal uniformly $J$-negative, see \cite[Corollary 1.5.2]{Ando}. Moreover, $Q\in \Q$.

Notice that $R(S(I-Q)^\#)=S(\M^{\ort})=S(S(\T)^{\ort})=S(S^{-1}(\T^{\ort}))=\T^{\ort}$. Therefore, $QS(I-Q)^\#=0$ and
\[
QS=QSQ^\# + QS(I-Q)^\#=QSQ^\#.
\]
Furthermore, if $\K{Sf}{f}\geq 0$ for every $f\in\T$ then $QS$ is $J$-positive. Analogously, if $\K{Sg}{g}\leq 0$ for every $g\in S(\T)^{\ort}$ then $(I-Q)S$ is $J$-negative. Then, by Theorem \ref{carac S}, $S$ is a $J$-frame operator.
  \end{pf*}

As it was proved in Proposition \ref{props S}, if an operator $S\in L(\HH)$ is a $J$-frame operator then it is an invertible $J$-selfadjoint operator satisfying $\ind_\pm(S)= \dim(\HH_\pm)$. Unfortunatelly, the converse is not true.

\begin{exmp}
Consider the Krein space obtained by endowing $\CC^2$ with the sesquilinear form 
\[
[(x_1,x_2),(y_1,y_2)]= x_1\overline{y_1} - x_2\overline{y_2},
\]
and the invertible $J$-selfadjoint operator $S$, whose matrix in the standard orthonormal basis is given by
\[
S=\matriz{0}{i}{i}{0}.
\]
Then, $S$ satisfies $\ind_\pm(S)= \dim(\HH_\pm)$, but it maps each $J$-positive vector into a $J$-negative vector. Then, by Proposition \ref{manda pos en pos}, $S$ cannot be a $J$-frame operator.
\end{exmp}

\section{Final remarks}

The following are some simple consequences of the material studied in the previous sections. Nevertheless, they are not going to be thoroughly developed in this notes.

\subsection*{Synthesis operators of $J$-frames as sums of plus and minus operators}

If $\F=\{f_i\}_{i\in I}$ is a $J$-frame for the Krein space $(\HH,\K{\,}{\,})$, it is easy to see that 
 $T_+$ and $T_+^\#$ are plus operators (considering $\ell_2(I)$ as a Krein space with the fundamental symmetry $J_2$ defined in \eqref{J2}). Furthermore, $T_+^\#$ is strict, and, $T_+$ is a strict plus operator if and only if $N(T)\cap \ell_2(I_+)=\{0\}$.

Also, these conditions have a natural counterpart for the operators $T_-$ and $T_-^\#$. Indeed, it follows analogously that $T_-$ and $T_-^\#$ are minus operators; $T_-^\#$ is always strict, and, $T_-$ is a strict minus operator if and only if $N(T)\cap \ell_2(I_-)=\{0\}$  (see \cite[Ch. 2]{Iokhvidov} for the terminology).

\subsection*{Frames for regular subspaces of a Krein space}

Given a Krein space $(\HH, \K{\,}{\,})$, recall that a subspace $\St$ of $\HH$ is \emph{regular} if there exists a (unique) $J$-selfadjoint projection onto $\St$. Since a regular subspace $\St$, endowed with the restriction of the indefinite inner product $\K{\,}{\,}$ to $\St$, is a Krein space (see \cite[Ch. 1,Theorem 7.16]{Iokhvidov}) the definition of $J$-frames applies for regular subspaces of $\HH$ too. Therefore, it is easy to infer a notion of ``$J$-frames for regular subspaces'' of a Krein space.

{\bf Juan I. Giribet}

jgiribet@fi.uba.ar

Departamento de Matem\'atica, FI-UBA, Buenos Aires, Argentina 

and 

IAM-CONICET.

\medskip

{\bf Alejandra Maestripieri}

amaestri@fi.uba.ar

Departamento de Matem\'atica, FI-UBA, Buenos Aires, Argentina 

and 

IAM-CONICET.

\medskip

{\bf Francisco Mart\'{\i}nez Per\'{\i}a}

francisco@mate.unlp.edu.ar

Departamento de Matem\'atica, FCE-UNLP, La Plata, Argentina 

and 

IAM-CONICET.

\medskip

{\bf Pedro G. Massey}

massey@mate.unlp.edu.ar

Departamento de Matem\'atica, FCE-UNLP, La Plata, Argentina 

and 

IAM-CONICET, Saavedra 15, Piso 3, (1083) Buenos Aires, Argentina.

\end{document}